\documentclass[12pt,a4paper]{article}
\usepackage{graphicx}
\usepackage{amsmath,amsfonts,amsxtra,amssymb,latexsym, amscd,amsthm}
\usepackage[mathscr]{eucal}
\usepackage{indentfirst}

\numberwithin{equation}{section}
\newtheorem{theorem}{Theorem}[section]
\newtheorem{lemma}[theorem]{Lemma}
\newtheorem{remark}[theorem]{Remark}

 3

\begin{document}
\title{The regularity and exponential decay of solution for a linear wave equation associated with two-point boundary conditions}
%THE REGULARITY AND EXPONENTIAL DECAY OF SOLUTION FOR  A  LINEAR WAVE EQUATION ASSOCIATED WITH TWO-POINT BOUNDARY CONDITIONS} %
\author{Le Xuan Truong $^{(1)}$, Le Thi Phuong Ngoc$^{(2)}$, \\
			Alain Pham Ngoc Dinh$^{(3)}$, Nguyen Thanh Long$^{(4)}$ }
\date{}
\maketitle
{\small
\begin{list}{}{
\setlength{\rightmargin}{.4in}
\setlength{\leftmargin}{.4in}
\setlength{\parsep}{.2pt}}
\item[$^{(1)}$] \em{University of Technical Education in HoChiMinh City, 01 Vo Van Ngan Str., Thu Duc Dist., HoChiMinh City, Vietnam. }\\
\hspace{0.2cm}{E-mail: lxuantruong@gmail.com }
\item[$^{(2)}$] Nhatrang Educational College, 01 Nguyen Chanh Str., Nhatrang City, Vietnam.\\
\hspace{0.2cm}{E-mail: ngoc1966@gmail.com}
\item[$^{(3)}$] MAPMO, UMR 6628, b\^{a}t. Math\'{e}matiques, University of Orl\'{e}ans, BP 6759, 45067 Orl\'{e}eans Cedex 2, France.\\
\hspace{0.2cm}{E-mail: alain.pham@univ-orleans.fr, alain.pham@math.cnrs.fr}
\item[$^{(4)}$] Department of Mathematics and Computer Science, University of Natural Science, Vietnam National University HoChiMinh City, 227 Nguyen Van Cu Str., Dist.5, HoChiMinh City, Vietnam.\\
\hspace{0.2cm}{E-mail: longnt@hcmc.netnam.vn, longnt2@gmail.com }
\end{list}
}
%\vspace{0.2cm}
\hrule
%\vspace{0.1cm}
{\small
{\center{\bf Abstract}}

%\vspace{0.3cm}
This paper is concerned with the existence and the regularity of global solutions to the linear wave equation associated the boundary conditions of two-point type. We also investigate the decay properties of the global solutions to this problem by the construction of a suitable Lyapunov functional.

\vspace{0.2cm}
%\vspace{0.1cm}

\vspace{0.3cm}
\hrule

\section{Introduction}
The wave equation
\begin{equation*}
	u_{tt}-\Delta u = f(x,t, u, u_{t}),
\end{equation*}
associated with the different boundary conditions, has been extensively studied by many authors, see [1 -8] and references therein. In the above mentioned papers, the existence and regularity of solutions, the asymptotic behavior and asymptotic expansion of solutions have received much attention. \\
\indent
In [8], Santos also studied the asymptotic behavior of the solutions to a coupled system of wave  equations having integral convolutions as memory terms. Their main result showed that the solution of that system decays uniformly in time, with rates depending on the rate of decay of the kernel
of the convolutions. \\
\indent
In this paper we consider the following initial-boundary value problem for the linear wave equation
\begin{equation}\label{eq1.01}
		u_{tt}-u_{xx}+Ku+\lambda u_{t}=f(x,t)\hspace{0.2cm} \text{in} \hspace{0.2cm} (0,1)\times (0,\infty),  \tag{1.1}
\end{equation}
\begin{equation}\label{eq1.02}
		u_{x}(0,t)=h_{0}u(0,t)+\lambda _{0}u_{t}(0,t)+\widetilde{h}_{1}u(1,t)+\widetilde{\lambda }_{1}u_{t}(1,t)+g_{0}(t)   \text{,}  \tag{1.2}
\end{equation}
\begin{equation}\label{eq1.03}
		-u_{x}(1,t)=h_{1}u(1,t)+\lambda _{1}u_{t}(1,t)+\widetilde{h}_{0}u(0,t)+\widetilde{\lambda }_{0}u_{t}(0,t)+g_{1}(t)  \text{,} \tag{1.3}
\end{equation}
\begin{equation}\label{eq1.04}
		u(x,0)=\widetilde{u}_{0}(x)\text{, \ }u_{t}(x,0)=\widetilde{u}_{1}(x), \tag{1.4}
\end{equation}
where $h_{0}$, $h_{1}$, $\lambda _{0}$, $\lambda _{1}$, $\widetilde{h}_{0}$, $\widetilde{h}_{1}$, $\widetilde{\lambda }_{0}$, $\widetilde{\lambda }_{1}$, $K$, $\lambda$ are constants and $\widetilde{u}_{0}$, $\widetilde{u}_{1}$, $f$, $g_{0}$,\ $g_{1}$ are given functions.

The rest of this paper consists of four sections. In section 2, we present some notations and lemmas that will be used to establish our results. In section 3, we investigate the existence and uniqueness of a weak solution and so a strong solutions of problem $\eqref{eq1.01} - \eqref{eq1.04}$ with the convenable conditions. Section 4 is devoted to the study of the regularity of solutions. Finally, in fifth section, we prove that the exponential decay properties of the global solutions are similar to that of the functionals $f,\ g_{0}, \ g_{1}$.
	
\section{Preliminaries}
	Let $\Omega = (0,1)$ and $Q_{T}=\Omega\times (0,T)$, for $T>0$. In what follows we will denote
\begin{equation*}
	\left< u, v\right> = \int_{0}^{1}u(x)v(x)dx, \quad \|v\| = \sqrt {\left< v, v\right>},
\end{equation*}
and $\|\cdot \|_{1}$ is a equivalent norm in $H^{1}(\Omega)$, defined by
\begin{equation*}
	\left\Vert v\right\Vert _{1}=\left( v^{2}(0)+ \|v_x\|^{2}  \right) ^{1/2}.
\end{equation*}
We also denote $u(x,t)$, $u_{t}(x,t)$, $u_{tt}(x,t)$, $u_{x}(x,t)$ and $u_{xx}(x,t)$ by $u(t)$, $u'(t)$, $u''(t)$, $u_{x}(t)$, $u_{xx}(t)$, respectively, when no confusion arises.

\vspace{0.2cm}
Consider a symmetric bilinear form $a(u,v)$ on $H^{1}(\Omega) \times H^{1}(\Omega)$ by setting
\begin{equation*}
	a(u,v)=\left<u_{x},v_{x}\right> + h_{0}u(0)v(0) + h_{1}u(1)v(1).
\end{equation*}
We state here some preliminary results that will be used in the sequel.
\begin{lemma}
	The imbedding $H^{1}(\Omega)\hookrightarrow C^{0}(\overline{\Omega })$ is compact and
\begin{equation}\label{eq2.01}
	\left\Vert v\right\Vert _{C^{0}(\overline{\Omega })}\leq \sqrt{2}\left\Vert v\right\Vert _{1}\text{, \ \textit{for all} }v\in H^{1}(\Omega). \tag{2.1}
\end{equation}

\end{lemma}
\begin{lemma}
	Let $h_{0}>0$ and $h_{1}\geq 0$. Then, the symmetric bilinear form $a(\cdot,\cdot)$ is continuous on $H^{1}(\Omega) \times H^{1}(\Omega)$ and coercive on $H^{1}(\Omega)$, i.e.,

\vspace{0.2cm}
\begin{minipage}[h]{0.8\textwidth}
\begin{description}
\item[(i)] $\left|a(u,v)\right| \leq C_{1}\|u\|_{1}\|v\|_{1}, \quad \forall u, v \in H^{1}(\Omega)$,
\item[(ii)] $a(v,v)\geq C_{0}\|v\|^2_{1}, \quad \forall  v \in H^{1}(\Omega)$,
\end{description}
\end{minipage}

\vspace{0.2cm} \noindent
where $C_{0}=\min\left\{1, h_{0}\right\}$ and $C_{1}=\max\left\{1,h_{0},2h_{1}\right\}$.
\end{lemma}

\begin{lemma}
	Let $\lambda _{0}$, $\lambda _{1}>0$ and $\widetilde{\lambda }_{0}$, $\widetilde{\lambda }_{1}\in \mathbb{R}$ such that $\left( \widetilde{\lambda }_{0}+\widetilde{\lambda }_{1}\right) ^{2}-4\lambda _{0}\lambda _{1}<0$. Then we have
\begin{equation*}
	\lambda _{0}x^{2}+\,\lambda _{1}y^{2}+\left( \widetilde{\lambda }_{0}+\widetilde{\lambda }_{1}\right) xy\geq \frac{1}{2}\mu _{\min }\left(x^{2}+y^{2}\right), \quad \forall x, y \in \mathbb{R},
\end{equation*}
where
\begin{equation*}
	\mu _{\min }=\frac{1}{4}\left[ -\left( \widetilde{\lambda }_{0}+\widetilde{\lambda }_{1}\right) ^{2}+4\,\lambda _{0}\lambda _{1}\right] \min \left\{ \frac{1}{\lambda _{0}}\text{, }\frac{1}{\lambda _{1}}\right\} >0.
\end{equation*}
\end{lemma}
\noindent The proof of these lemmas are straightforward. We shall omit the details.

\begin{remark}
	From the Lemma 2.2 we deduce that
\begin{equation}\label{eq2.02}
	C_{0}\|v\|^{2}_{1} \leq \|v\|^{2}_{a} \leq C_{1}\|v\|^{2}_{1}, \quad \forall v\in H^{1}(\Omega), \tag{2.2}
\end{equation}
where $\|\cdot\|_{a}$ is the norm on $H^{1}(\Omega)$ generated by the the symmetric bilinear form $a(\cdot,\cdot)$, i.e.,
\begin{equation*}
	\|v\|_{a}=\sqrt{a(v,v)}, \hspace{0.2cm} \forall v\in H^{1}(\Omega).
\end{equation*}
\end{remark}

\section{Existence and uniqueness of solutions}
In this section, we assume that $h_{0}$, $\lambda_{0}$, $\lambda_{1}$ are positive constants, $h_{1}$ is nonnegative constant and $K$, $\lambda$, $\widetilde{h}_{0}$, $\widetilde{h}_{1}$, $\widetilde{\lambda }_{0}$, $\widetilde{\lambda }_{1}$ are constants verifying the condition
\begin{equation}\label{eq3.01}
	\left|\widetilde{\lambda }_{0}+\widetilde{\lambda }_{1}\right|< 2\sqrt{\lambda_{0}\lambda_{1}}. \tag{3.1}
\end{equation}

\begin{theorem}
	Let $T>0$ and assume that $g_{0}\text{, }g_{1}\in L^{2}(0,T)$, $f \in L^{1}\left( 0,T;L^{2}(\Omega)\right)$. Then, for each $\left(\widetilde{u}_{0}, \widetilde{u}_{1}\right) \in H^{1}(\Omega) \times L^{2}(\Omega)$, the problem $\eqref{eq1.01} - \eqref{eq1.04}$ has a unique weak solution $u$ satisfying
\begin{equation*}
	u\in L^{\infty }\left(0,T;H^{1}(\Omega)\right), \hspace{0.2cm} u_{t}\in L^{\infty }\left(0,T;L^{2}(\Omega)\right),
\end{equation*}
and
\begin{equation*}
	u(0,\cdot ), u(1,\cdot )\in H^{1}\left( 0,T\right).
\end{equation*}
\end{theorem}

\begin{proof}
The proof consists of step 1 - 4.

\vspace{0.2cm}
\noindent {\bf Step 1.} The Faedo-Galerkin approximation. Let $\{w_{j}\}$\ be a denumerable base of $\ H^{1}(\Omega)$. We find the approximate solution of problem $\eqref{eq1.01} - \eqref{eq1.04}$ in the form
\begin{equation*}
	u_{m}(t)=\sum_{j=1}^{m}c_{mj}(t)w_{j},
\end{equation*}
where the coefficient functions $c_{mj}$\ satisfy the system of ordinary differential equations
\begin{gather}
	\left< u''_{m}(t),w_{j}\right>+a\left(u_{m}(t),w_{j}\right) +\left( \lambda _{0}u'_{m}(0,t)+\widetilde{h}%
_{1}u_{m}(1,t)+\widetilde{\lambda }_{1}u'_{m}(1,t)\right) w_{j}(0)\notag\\
	+ \left( \lambda _{1}u'_{m}(1,t)+\widetilde{h}_{0}u_{m}(0,t)+\widetilde{\lambda }_{0}u'_{m}(0,t)\right) w_{j}(1)+\left< Ku_{m}(t)+\lambda u'_{m}(t),w_{j}\right> \label{eq3.02} \tag{3.2} \\
	=-g_{0}(t)w_{j}(0)-g_{1}(t)w_{j}(1)+\left<f(t),w_{j}\right>,  1\leq j \leq m, \notag
\end{gather}
with the initial conditions
\begin{equation}\label{eq3.03}
	u_{m}(0)=u_{0m}=\sum_{j=1}^{m}\alpha _{mj}w_{j}\rightarrow \widetilde{u}_{0} \hspace{0.2cm} \text{strongly in} \hspace{0.2cm} H^1(\Omega), \tag{3.3}
\end{equation}
and
\begin{equation}\label{eq3.04}
	u'_{m}(0)=u_{1m}=\sum_{j=1}^{m}\beta_{mj}w_{j}\rightarrow \widetilde{u}_{1} \hspace{0.2cm} \text{strongly in} \hspace{0.2cm} L^2(\Omega). \tag{3.4}
\end{equation}
From the assumptions of Theorem 3.1, system $\eqref{eq3.02}-\eqref{eq3.04}$ has solution $u_{m}(t)$ on some interval $[0,T_{m}]$. The following estimates allow one to take $T_{m}=T$, for all $m$.

\vspace{0.2cm}
\noindent {\bf Step 2.} A priori estimates. Multiplying the $j^{th}$ equation of $\eqref{eq3.02}$ by $c_{mj}^{\prime}(t)$ and summing up with respect to $j$, afterwards, integrating by parts with respect to the time variable from $0$ to $t$, we
get after some rearrangements
\begin{alignat}{2}\label{eq3.05}
	S_{m}(t)=&S_{m}(0)-2\widetilde{h}_{0}\int_{0}^{t}u_{m}(0,s)u'_{m}(1,s)ds -2\widetilde{h}_{1}\int_{0}^{t}u_{m}(1,s)u'_{m}(0,s)ds \notag \\
	& -2K\int_{0}^{t}\left\langle u_{m}(s),u'_{m}(s)\right\rangle ds-2\lambda \int_{0}^{t}\left\Vert u'_{m}(s)\right\Vert ^{2}ds \tag{3.5} \\
	& -2\int_{0}^{t}g_{0}(s)u'_{m}(0,s)ds-2\int_{0}^{t}g_{1}(s)u'_{m}(1,s)ds+2\int_{0}^{t}\left\langle f(s),u'_{m}(s)\right\rangle ds \notag\\
	=&S_{m}(0)+\sum_{i=1}^{7}I_{i}, \notag
\end{alignat}
where
\begin{alignat}{2}\label{eq3.06}
	S_{m}(t)=& \left\Vert u'_{m}(t)\right\Vert ^{2}+\left\Vert u_{m}(t)\right\Vert _{a}^{2} \notag \\
	&+2\int_{0}^{t}\left[ \lambda _{0}\left\vert u'_{m}(0,s)\right\vert^{2}+\lambda _{1}\,\left\vert u'_{m}(1,s)\right\vert ^{2}+\left( \widetilde{\lambda }_{0}+\widetilde{\lambda }_{1}\right)u'_{m}(0,t)u'_{m}(1,t)\right] ds. \tag{3.6}
\end{alignat}
By Lemma 2.3, it follows from $\eqref{eq3.06}$, that
\begin{equation}\label{eq3.07}
	S_{m}(t) \geq \mu_{0}X_{m}(t), \tag{3.7}
\end{equation}
where
\begin{equation}\label{eq3.08}
	X_m(t) = \left\Vert u'_{m}(t)\right\Vert ^{2}+\left\Vert u_{m}(t)\right\Vert _{1}^{2}+\int_{0}^{t}\left( \left\vert
u'_{m}(0,s)\right\vert ^{2}+\left\vert u'_{m}(1,s)\right\vert^{2}\right) ds, \tag{3.8}
\end{equation}
and $\mu_{0} = \min\left\{C_{0}, \mu_{min}\right\}$.

\vspace{0.2cm}
\noindent Now, using the inequalities $\eqref{eq2.01}-\eqref{eq2.02}$ and the following inequalities
\begin{equation}\label{eq3.09}
	2ab\leq \varepsilon a^{2}+\frac{1}{\varepsilon }b^{2}, \forall a,b\in \mathbb{R}, \forall \varepsilon >0, \tag{3.9}
\end{equation}
\begin{equation}\label{eq3.10}
	\left\vert u_{m}(0,t)\right\vert \leq \left\Vert u_{m}(t)\right\Vert_{C^{0}(\overline{\Omega })}\leq \sqrt{2}\left\Vert u_{m}(t)\right\Vert_{1}\leq \sqrt{2X_{m}(t)}, \tag{3.10}
\end{equation}
\begin{equation}\label{eq3.11}
	\left\Vert u_{m}(t)\right\Vert ^{2}\leq 2\left\Vert u_{0m}\right\Vert^{2}+2\int_{0}^{t}\left\Vert u'_{m}(s)\right\Vert ^{2}ds\leq 2\left\Vert u_{0m}\right\Vert ^{2}+2\int_{0}^{t}X_{m}(s)ds, \tag{3.11}
\end{equation}
we shall estimate respectively the terms on the right-hand side of $\eqref{eq3.05}$ as follows
\begin{alignat}{2}\label{eq3.12}
	I_{1}&=-2\widetilde{h}_{0}\int_{0}^{t}u_{m}(0,s)u'_{m}(1,s)ds \notag\\
	&\leq \frac{1}{\varepsilon }\left\vert \widetilde{h}_{0}\right\vert^{2}\,\int_{0}^{t}\left\vert u_{m}(0,s)\right\vert ^{2}ds+\varepsilon\int_{0}^{t}\left\vert u'_{m}(1,s)\right\vert ^{2}ds \tag{3.12} \\
	& \leq \frac{2}{\varepsilon }\left\vert \widetilde{h}_{0}\right\vert ^{2}\int_{0}^{t}X_{m}(s)ds+\varepsilon X_{m}(t), \notag
\end{alignat}
\begin{equation}\label{eq3.13}
	I_{2}=-2\widetilde{h}_{1}\int_{0}^{t}u_{m}(1,s)u'_{m}(0,s)ds\leq \frac{2}{\varepsilon }\left\vert \widetilde{h}_{1}\right\vert ^{2}\int_{0}^{t}X_{m}(s)ds+\varepsilon X_{m}(t), \tag{3.13}
\end{equation}
\begin{equation}\label{eq3.14}
	I_{3}=-2K\int_{0}^{t}\left\langle u_{m}(s),u'_{m}(s)\right\rangle ds\leq2\sqrt{2}\left\vert K\right\vert \int_{0}^{t}X_{m}(s)ds, \tag{3.14}
\end{equation}
\begin{equation}\label{eq3.15}
	I_{4}=-2\lambda \int_{0}^{t}\left\Vert u'_{m}(s)\right\Vert ^{2}ds\leq 2\left\vert \lambda \right\vert \int_{0}^{t}X_{m}(s)ds, \tag{3.15}
\end{equation}
\begin{alignat}{2}\label{eq3.16}
	I_{5}=-2\int_{0}^{t}g_{0}(s)u'_{m}(0,s)ds &\leq \frac{1}{\varepsilon }\left\Vert g_{0}\right\Vert _{L^{2}\left( 0,T\right) }^{2}+\varepsilon\int_{0}^{t}\left\vert u'_{m}(0,s)\right\vert ^{2}ds \notag \\
	&\leq \frac{1}{\varepsilon }\left\Vert g_{0}\right\Vert _{L^{2}\left(0,T\right) }^{2}+\varepsilon X_{m}(t), \tag{3.16}
\end{alignat}
\begin{equation}\label{eq3.17}
	I_{6}=-2\int_{0}^{t}g_{1}(s)u'_{m}(1,s)ds\leq \frac{1}{\varepsilon }\left\Vert g_{1}\right\Vert _{L^{2}\left( 0,T\right) }^{2}+\varepsilon X_{m}(t), \tag{3.17}
\end{equation}
\begin{equation}\label{eq3.18}
	I_{7}=2\int_{0}^{t}\left\langle f(s),u'_{m}(s)\right\rangle ds\leq \int_{0}^{T}\left\Vert f(s)\right\Vert ds+\int_{0}^{t}\left\Vert f(s)\right\Vert X_{m}(s)ds. \tag{3.18}
\end{equation}
On the other hand, using $\eqref{eq3.03} -\eqref{eq3.04}$, $\eqref{eq3.06}$ and the assumption $\left(\widetilde{u}_{0}, \widetilde{u}_{1}\right) \in H^{1}(\Omega) \times L^{2}(\Omega)$, we have
\begin{equation}\label{eq3.19}
	S_{m}(0)=\left\Vert u_{1m}\right\Vert ^{2}+\left\Vert u_{0m}\right\Vert_{a}^{2}\leq \widetilde{C}_{1}\text{ for all\textrm{\ }}m, \tag{3.19}
\end{equation}
where $\widetilde{C}_{1}$ is a constant depending only on $\widetilde{u}_{0} $, $\widetilde{u}_{1}$, $h_{0}$ and $h_{1}$.

\vspace{0.2cm} Combining $\eqref{eq3.05}$, $\eqref{eq3.07}$, $\eqref{eq3.12}$-$\eqref{eq3.19}$, we obtain
\begin{alignat}{2}\label{eq3.20}
	\left(\mu_{0}-4\varepsilon\right)X_{m}(t) &\leq \widetilde{C}_{1}+\frac{1}{\varepsilon }\left\Vert g_{0}\right\Vert _{L^{2}\left(0,T\right) }^{2}+\frac{1}{\varepsilon }\left\Vert g_{1}\right\Vert _{L^{2}\left( 0,T\right) }^{2}+\int_{0}^{T}\left\Vert f(s)\right\Vert ds \tag{3.20}\\
	&+\int_{0}^{t}\left[\frac{2}{\varepsilon}\left(\left\vert \widetilde{h}_{0}\right\vert ^{2}+\left\vert \widetilde{h}_{1}\right\vert ^{2}\right)+2\sqrt{2}\left\vert K\right\vert+2\left\vert \lambda \right\vert+ \|f(s)\|\right]X_{m}(s)ds, \notag
\end{alignat}
for all $\varepsilon >0$. By choosing $\varepsilon >0$ such that $\mu_{0}-4\varepsilon >0$, it follows from $\eqref{eq3.20}$ that
\begin{equation}\label{eq3.21}
	X_{m}(t)\leq M_{T}^{(1)}+\int_{0}^{t}N_{T}^{(1)}(s)X_{m}(s)ds, \tag{3.21}
\end{equation}
where
\begin{equation*}
	M_{T}^{(1)}=\left( \mu _{0}-4\varepsilon \right) ^{-1}\left( \widetilde{C}_{1}+\frac{1}{\varepsilon }\left\Vert g_{0}\right\Vert _{L^{2}\left(0,T\right) }^{2}+\frac{1}{\varepsilon }\left\Vert g_{1}\right\Vert_{L^{2}\left( 0,T\right) }^{2}+\int_{0}^{T}\left\Vert f(s)\right\Vert ds\right),
\end{equation*}
and
\begin{equation*}
	N_{T}^{(1)}(s)=\left( \mu _{0}-4\varepsilon \right) ^{-1}\left[ \frac{2}{\varepsilon }\left( \left\vert \widetilde{h}_{0}\right\vert^{2}+\left\vert \widetilde{h}_{1}\right\vert ^{2}\right) +2\sqrt{2}\left\vert K\right\vert +2\left\vert \lambda \right\vert+\left\Vert f(s)\right\Vert \right], \ N_{T}^{(1)}\in L^{1}(0,T).
\end{equation*}
By Gronwall's lemma, we deduce from $\eqref{eq3.21}$, that
\begin{equation}\label{eq3.22}
	X_{m}(t)\leq M_{T}^{(1)}\exp \left( \int_{0}^{t}N_{T}^{(1)}(s)ds\right) \leq C_{T}\text{, \ for all }t\in [0,T], \tag{3.22}
\end{equation}
where $C_{T}$ is a posistive constant depending only on $T$.

\vspace{0.2cm}
\noindent {\bf Step 3.} Limiting process. From $\eqref{eq3.08}$ and $\eqref{eq3.22}$, we deduce the existence of a subsequence of $\{u_{m}\}$\ still also so denoted, such that
\begin{equation}\label{eq3.23}
	\left\{
\begin{array}{cccc}
u_{m}\rightarrow u & \text{in} & L^{\infty }(0,T;H^{1}(\Omega))\text{ } & \text{%
weak*,} \\
u'_{m}\rightarrow u' & \text{in} & L^{\infty }(0,T;L^{2}(\Omega)) & \text{%
weak*,} \\
u_{m}(0,\cdot )\rightarrow u(0,\cdot ) & \text{in} & H^{1}(0,T) & \text{%
weakly,} \\
u_{m}(1,\cdot )\rightarrow u(1,\cdot ) & \text{in} & H^{1}(0,T) & \text{%
weakly.}%
\end{array}%
\right. \tag{3.23}
\end{equation}
By the compactness lemma of \ Lions [5: p.57] and the imbedding $H^{1}(0,T)$  $\hookrightarrow C^{0}\left( \left[ 0,T\right] \right) $, we can deduce from $\eqref{eq3.23}$ the existence of a subsequence still denoted by $\{u_{m}\}$, such that
\begin{equation}\label{eq3.24}
	\left\{
\begin{array}{ccc}
u_{m}\rightarrow u & \text{strongly\thinspace \thinspace in} & L^{2}(Q_{T})%
\text{,} \\
u_{m}(0,\cdot )\rightarrow u(0,\cdot ) & \text{strongly\thinspace \thinspace
in} & C^{0}\left( [0,T]\right) \text{,} \\
u_{m}(1,\cdot )\rightarrow u(1,\cdot ) & \text{strongly\thinspace \thinspace
in} & C^{0}\left( [0,T]\right) \text{.}%
\end{array}%
\right. \tag{3.24}
\end{equation}
Passing to the limit in $\eqref{eq3.02}$-$\eqref{eq3.04}$ by $\eqref{eq3.23}$ and $\eqref{eq3.24}$ we have $u$
satisfying the equation
\begin{alignat}{2}\label{eq3.25}
	\frac{d}{dt}\left\langle u'(t),v\right\rangle &+a\left( u(t),v\right)+\left( \lambda _{0}u'(0,t)+\widetilde{h}_{1}u(1,t)+\widetilde{\lambda }_{1}u'(1,t)\right) v(0) \notag \\
	&+\left( \lambda _{1}u'(1,t)+\widetilde{h}_{0}u(0,t)+\widetilde{\lambda }_{0}u'(0,t)\right) v(1)+\left\langle Ku(t)+\lambda u'(t),v\right\rangle \tag{3.25}\\
	&= -g_{0}(t)v(0)-g_{1}(t)v(1)+\left\langle f(t), v\right\rangle \text{, \ }\forall \,\,v\in H^{1}(\Omega), \notag
\end{alignat}
in $L^{2}(0,T)$ weakly, and
\begin{equation}\label{eq3.26}
	u(0)=\widetilde{u}_{0}\text{, }\ u'(0)=\widetilde{u}_{1}. \tag{3.26}
\end{equation}
The existence of the theorem is proved completely.

\vspace{0.2cm}
\noindent {\bf Step 4.} Uniqueness of the solution. Let $u_{1}$, $u_{2}$ be two weak solutions of problem $\eqref{eq1.01} - \eqref{eq1.04}$, such that
\begin{equation}\label{eq3.27}
	\left\{
	\begin{array}{l}
		u_{i}\in L^{\infty }(0,T;H^{1}(\Omega))\text{, }u'_{i}\in L^{\infty }(0,T;L^{2}(\Omega)), \\
		u_{i}(0,\cdot )\text{, }u_{i}(1,\cdot )\in \, H^{1}\left( 0,T\right), i =1, 2.
	\end{array}
	\right. \tag{3.27}
\end{equation}
Then $u=u_{1}-u_{2}$ is the weak solution of the following problem
\begin{equation}\label{eq3.28}
	\left\{
	\begin{array}{l}
		u_{tt}-u_{xx}+Ku+\lambda u_{t}=0\text{, \ }(x, t)\in Q_{T}, \\
		u_{x}(0,t)=h_{0}u(0,t)+\lambda _{0}u_{t}(0,t)+\widetilde{h}_{1}u(1,t)+\widetilde{\lambda }_{1}u_{t}(1,t),\\
		-u_{x}(1,t)=h_{1}u(1,t)+\lambda _{1}u_{t}(1,t)+\widetilde{h}_{0}u(0,t)+\widetilde{\lambda }_{0}u_{t}(0,t),\\
		u(x,0)=0\text{, \ }u_{t}(x,0)=0.
	\end{array}
	\right. \tag{3.28}
\end{equation}
By using the lemma in [8, Lemma 2.4, p. 1799], we deduce that
\begin{alignat}{2}\label{eq3.29}
	\left\Vert u'(t)\right\Vert ^{2}&+\left\Vert u(t)\right\Vert _{a}^{2}+2\int_{0}^{t}\left\langle Ku(s)+\lambda u'(s),u'(s)\right\rangle ds \notag\\
	&+2\int_{0}^{t}\left[ \lambda _{0}\left\vert u'(0,s)\right\vert^{2}ds+\lambda _{1}\,\left\vert u'(1,s)\right\vert ^{2}ds+\left( \widetilde{\lambda }_{0}+\widetilde{\lambda }_{1}\right) u'(0,s)u'(1,s)\right] ds \tag{3.29}\\
	&+2\widetilde{h}_{1}\int_{0}^{t}u(1,s)u'(0,s)ds+2\widetilde{h}_{0}\int_{0}^{t}u(0,s)u'(1,s)ds. \notag
\end{alignat}
Putting
\begin{equation}\label{eq3.30}
	\sigma (t)=\left\Vert u'(t)\right\Vert ^{2}+\left\Vert u(t)\right\Vert _{1}^{2}+\mu _{\min }\,\int_{0}^{t}\left[ \left\vert u'(0,s)\right\vert ^{2}+\left\vert u'(1,s)\right\vert ^{2}\right] ds. \tag{3.30}
\end{equation}
From $\eqref{eq3.29}$, $\eqref{eq3.30}$ and Lemma 2.3, we prove, in a similar manner to that in the above part, that
\begin{equation}\label{eq3.31}
	\left( 1-\frac{2\varepsilon }{\mu _{\min }}\right) \sigma (t)\leq 2\left[ \frac{1}{\varepsilon }\left( \left\vert \widetilde{h}_{0}\right\vert^{2}+\left\vert \widetilde{h}_{1}\right\vert ^{2}\right) +\sqrt{2}\left\vert K\right\vert +\left\vert \lambda \right\vert \right]\int_{0}^{t}\sigma (s)ds. \tag{3.31}
\end{equation}
Choosing $\varepsilon >0$, with $1-2\varepsilon \mu _{\min }^{-1}>0$. Using Gronwall's lemma, it follows from $\eqref{eq3.30}$-$\eqref{eq3.31}$, that $\sigma (t)\equiv 0$, i.e., $u_{1}\equiv u_{2}$. The theorem 3.1 is proved
completely.
\end{proof}

\begin{theorem}
	Let $T>0$ and assume that $g_{0}\text{, }g_{1}\in H^{1}(0,T)$, $f , f_{t}\in L^{2}\left( Q_{T}\right)$. Then, for each $\left(\widetilde{u}_{0}, \widetilde{u}_{1}\right) \in H^{2}(\Omega) \times H^{1}(\Omega)$, the problem $\eqref{eq1.01} - \eqref{eq1.04}$ has a unique weak solution $u$ satisfying
\begin{equation}\label{eq3.32}
	u\in L^{\infty }\left(0,T;H^{2}(\Omega)\right), \hspace{0.2cm} u_{t}\in L^{\infty }\left(0,T;H^{1}(\Omega)\right), \hspace{0.2cm} u_{tt}\in L^{\infty }\left(0,T;L^{2}(\Omega)\right), \tag{3.32}
\end{equation}
and
\begin{equation}\label{eq3.33}
	u(0,\cdot ), u(1,\cdot )\in H^{2}\left( 0,T\right). \tag{3.33}
\end{equation}
\end{theorem}

\begin{proof}
	The proof consists of Steps 1-4.

\vspace{0.2cm}
\noindent {\bf Step 1.} The Faedo-Galerkin approximation. Let $\{w_{j}\}$ be a denumerable base of $H^{2}(\Omega).$ We find the approximate solution of problem $\eqref{eq1.01} - \eqref{eq1.04}$ in the form
\begin{equation*}
	u_{m}(t)=\sum_{j=1}^{m}c_{mj}(t)w_{j},
\end{equation*}
where the coefficient functions $c_{mj}$\ satisfy the system of ordinary differential equations $\eqref{eq3.02}$, with the initial conditions
\begin{equation}\label{eq3.34}
	\left\{
	\begin{array}{l}
	u_{m}(0)=u_{0m}=\sum_{j=1}^{m}\alpha _{mj}w_{j}\rightarrow \widetilde{u}_{0} \text{ strongly\thinspace \thinspace \thinspace \thinspace in\thinspace }\,\,\,H^{2}\text{,} \\
	u_{m}^{/}(0)=u_{1m}=\sum_{j=1}^{m}\beta _{mj}w_{j}\rightarrow \widetilde{u}_{1}\,\,\text{strongly\thinspace \thinspace \thinspace in\thinspace }\,\,\,H^{1}\text{. }
	\end{array}
	\right. \tag{3.34}
\end{equation}
From the assumptions of Theorem 3.2, system $\eqref{eq3.02}$ and $\eqref{eq3.34}$ has solution $u_{m}(t)$ on some interval $[0,T_{m}]$. The following estimates
allow one to take $T_{m}=T$\ \ for all $m$.

\vspace{0.2cm}
\noindent {\bf Step 2.} A priori estimates. By same arguments as in proof of Theorem 3.1, we obtain
\begin{equation}\label{eq3.35}
	X_{m}(t) \leq C_{T}, \, \text{for all} \, \, \, t\in [0,T], \, \, m \in \mathbb{Z}^{+}, \tag{3.35}
\end{equation}
where $X_{m}(t)$ defined by $\eqref{eq3.08}$ and $C_{T}$ always indicating a bound depending on $T$.

\vspace{0.2cm} Now, differentiating $\eqref{eq3.02}$ with respect to $t$, we have
\begin{gather}
	\left\langle u'''_{m}(t),w_{j}\right\rangle +a\left(u'_{m}(t),w_{j}\right) +\left( \lambda _{0}u''_{m}(0,t)+\widetilde{h}_{1}u'_{m}(1,t)+\widetilde{\lambda }_{1}u''_{m}(1,t)\right) w_{j}(0)+ \notag \\
	\left( \lambda _{1}u''_{m}(1,t)+\widetilde{h}_{0}u'_{m}(0,t)+\widetilde{\lambda }_{0}u''_{m}(0,t)\right) w_{j}(1)+\left\langle Ku'_{m}(t)+\lambda u''_{m}(t),w_{j}\right\rangle \label{eq3.36} \tag{3.36}\\
	=-g'_{0}(t)w_{j}(0)-g'_{1}(t)w_{j}(1)+\left\langle f^{\prime}(t), w_{j}~\right\rangle \text{,} \notag
\end{gather}
for all $j = 1, 2, ...,m$.

\vspace{0.2cm}
Multiplying the $j^{th}$ equation of $\eqref{eq3.36}$ by $c'_{mj}(t)$, summing up with respect to $j$ and then integrating with respect to the time variable from $0$\ to $t$, we have after some rearrangements
\begin{alignat}{2}\label{eq3.37}
	\widetilde{S}_{m}(t)=&\widetilde{S}_{m}(0)-2\widetilde{h}_{0} \int_{0}^{t}u'_{m}(0,s)u''_{m}(1,s)ds-2\widetilde{h}_{1}\,\int_{0}^{t}u'_{m}(1,s)u''_{m}(0,s)ds \notag \\
	&-2K\int_{0}^{t}\left\langle u'_{m}(s),u''_{m}(s)\right\rangle ds-2\lambda \int_{0}^{t}\left\Vert u''_{m}(s)\right\Vert ^{2}ds \tag{3.37}\\
	&-2\int_{0}^{t}g'_{0}(s)u''_{m}(0,s)ds-2\int_{0}^{t}g'_{1}(s)u''_{m}(1,s)ds+2\int_{0}^{t}\left\langle f^{\prime}(s),u''_{m}(s)\right\rangle ds \notag \\
	&=\widetilde{S}_{m}(0)+\sum_{i=1}^{7}J_{i}\text{, } \notag
\end{alignat}
where
\begin{alignat}{2}\label{eq3.38}
	\widetilde{S}_{m}(t)=&\left\Vert u''_{m}(t)\right\Vert ^{2}+\left\Vert u'_{m}(t)\right\Vert _{a}^{2} \tag{3.38}\\
	&+2\int_{0}^{t}\left[ \lambda _{0}\left\vert u''_{m}(0,s)\right\vert^{2}+\lambda _{1}\,\left\vert u''_{m}(1,s)\right\vert ^{2}+\left( \widetilde{\lambda }_{0}+\widetilde{\lambda }_{1}\right)\,u''_{m}(0,t)u''_{m}(1,t)\right] ds\text{.} \notag
\end{alignat}
Using $\eqref{eq3.34}$, $\eqref{eq3.38}$ and Lemma 2.1, we have
\begin{equation}\label{eq3.39}
	\widetilde{S}_{m}(0)=\left\Vert u''_{m}(0)\right\Vert ^{2}+\left\Vert u_{1m}\right\Vert _{a}^{2}\leq \widetilde{C}_{2}, \, \, \text{for all} \,\, m, \tag{3.39}
\end{equation}
where $\widetilde{C}_{2}$ is a constant depending only on $\widetilde{u}_{0} $, $\widetilde{u}_{1}$, $f(\cdot ,0)$, $K$ and $\lambda $. On the other hand, by Lemma 2.3, it follows from $\eqref{eq3.39}$ that
\begin{equation}\label{eq3.40}
	\widetilde{S}_{m}(t)\geq \mu _{0}\widetilde{X}_{m}(t)\text{,} \tag{3.40}
\end{equation}
where
\begin{equation}\label{eq3.41}
	\widetilde{X}_{m}(t)=\left\Vert u''_{m}(t)\right\Vert ^{2}+\left\Vert u'_{m}(t)\right\Vert _{1}^{2}+\int_{0}^{t}\left( \left\vert u''_{m}(0,s)\right\vert ^{2}+\left\vert u''_{m}(1,s)\right\vert^{2}\right) ds\text{,} \tag{3.41}
\end{equation}
and $ \mu _{0}=\min \{C_{0},\mu _{\min }\}$.

\vspace{0.2cm}
By estimating the terms $J_{i}, \, (i = 1, 2,..., 7)$ on the right-hand side of $\eqref{eq3.37}$ as in the proof of Theorem 3.1, we get
\begin{equation}\label{eq3.42}
	\widetilde{X}_{m}(t)\leq M_{T}^{(2)}+\int_{0}^{t}N_{T}^{(2)}(s)\widetilde{X}_{m}(s)ds, \tag{3.42}
\end{equation}
where
\begin{equation*}
	M_{T}^{(2)}=\frac{2}{\mu _{0}}\left[ \widetilde{C}_{2}+\frac{8}{\mu _{0}}\left\Vert g'_{0}\right\Vert _{L^{2}\left( 0,T\right) }^{2}+\frac{8}{\mu_{0}}\left\Vert g'_{1}\right\Vert _{L^{2}\left( 0,T\right)}^{2}+\int_{0}^{T}\left\Vert f^{\prime}(s)\right\Vert ds\right] \text{, }
\end{equation*}
and
\begin{equation*}
	N_{T}^{(2)}(s)=\frac{2}{\mu _{0}}\left[ \frac{16}{\mu _{0}}\left(\widetilde{h}_{0} ^{2}+\widetilde{h}_{1} ^{2}\right) +2\sqrt{2}\left\vert K\right\vert +2\left\vert \lambda \right\vert +\left\Vert f^{\prime}(s)\right\Vert \right],\  N_{T}^{(2)} \in L^{1}(0,T).
\end{equation*}
From $\eqref{eq3.42}$ and applying Gronwall's inequality, we obtain that
\begin{equation}\label{eq3.43}
	\widetilde{X}_{m}(t)\leq M_{T}^{(2)}\exp \left(\int_{0}^{t}N_{T}^{(2)}(s)ds\right) \leq C_{T}\text{, \ for all }t\in\lbrack 0,T]\text{.} \tag{3.43}
\end{equation}

\vspace{0.2cm}
\noindent {\bf Step 3.} Limiting process. From $\eqref{eq3.08}$, $\eqref{eq3.35}$, $\eqref{eq3.41}$ and $\eqref{eq3.43}$, we deduce the existence of a subsequence of $\{u_{m}\}$, still denoted by $\{u_{m}\}$, such that
\begin{equation}\label{eq3.44}
	\left\{
\begin{array}{llll}
u_{m}\rightarrow u & \text{in} & L^{\infty }(0,T;H^{1}(\Omega))\text{ } & \text{%
weak*,} \\
u'_{m}\rightarrow u' & \text{in} & L^{\infty }(0,T;H^{1}(\Omega)) & \text{%
weak*,} \\
u''_{m}\rightarrow u'' & \text{in} & L^{\infty }(0,T;L^{2}(\Omega)) & \text{%
weak*,} \\
u_{m}(0,\cdot )\rightarrow u(0,\cdot ) & \text{in} & H^{2}(0,T) & \text{%
weakly,} \\
u_{m}(1,\cdot )\rightarrow u(1,\cdot ) & \text{in} & H^{2}(0,T) & \text{%
weakly.}%
\end{array}%
\right. \tag{3.44}
\end{equation}
By the compactness lemma of Lions [5, p.57] and the imbeddings $H^{1}(0,T)$ $\hookrightarrow C^{0}\left( \left[ 0,T\right] \right) $, $H^{2}(0,T)$ $\hookrightarrow C^{1}\left( \left[ 0,T\right] \right) $, we can deduce from $\eqref{eq3.44}$ the existence of a subsequence still denoted by $\{u_{m}\}$, such that
\begin{equation}\label{eq3.45}
	\left\{
\begin{array}{lll}
u_{m}\rightarrow u & \text{strongly\thinspace \thinspace in} & L^{2}(Q_{T})%
\text{, \ and a.e.\ }(x,t)\in Q_{T}\text{,} \\
u_{m}^{/}\rightarrow u^{/} & \text{strongly\thinspace \thinspace in} &
L^{2}(Q_{T})\text{, \ and a.e.\ }(x,t)\in Q_{T}\text{,} \\
u_{m}(0,\cdot )\rightarrow u(0,\cdot ) & \text{strongly\thinspace \thinspace
in} & C^{1}\left( [0,T]\right) \text{,} \\
u_{m}(1,\cdot )\rightarrow u(1,\cdot ) & \text{strongly\thinspace \thinspace
in} & C^{1}\left( [0,T]\right) \text{.}%
\end{array}%
\right. \tag{3.45}
\end{equation}
Passing to the limit in $\eqref{eq3.02}$ and $\eqref{eq3.34}$ by $\eqref{eq3.44}$-$\eqref{eq3.45}$ we have $u$ satisfying the problem
\begin{gather}
	\left\langle u''(t),v\right\rangle +a\left( u(t),v\right) +\left(\lambda _{0}u'(0,t)+\widetilde{h}_{1}u(1,t)+\widetilde{\lambda }_{1}u'(1,t)\right) v(0) \notag\\
	+\left( \lambda _{1}u'(1,t)+\widetilde{h}_{0}u(0,t)+\widetilde{\lambda }_{0}u'(0,t)\right) v(1)+\left\langle Ku(t)+\lambda u'(t),v\right\rangle \label{eq3.46} \tag{3.46}\\
	=-g_{0}(t)v(0)-g_{1}(t)v(1)+\left\langle f\,(t),v\right\rangle \text{, \ \ }\forall \,\,v\in H^{1}(\Omega)\text{,} \notag
\end{gather}
\begin{equation}\label{eq3.47}
	u(0)=\widetilde{u}_{0}\text{, }\ u'(0)=\widetilde{u}_{1}\text{.} \tag{3.47}
\end{equation}
\indent On the other hand, it follows from $\eqref{eq3.44}_{1,2,3}$ and $\eqref{eq3.46}$, that
\begin{equation}\label{eq3.48}
	u_{xx}=u''+Ku+\lambda u'-f~\in L^{\infty }(0,T;L^{2}(\Omega))\text{.} \tag{3.48}
\end{equation}
Thus, $u\in L^{\infty }(0,T;H^{2}(\Omega))$\ and the existence of solution is proved completely.

\vspace{0.2cm}
\noindent {\bf Step 4.} Uniqueness of the solution of problem $\eqref{eq1.01} - \eqref{eq1.04}$ is similarly proved as in  Theorem 3.1 and we will omit here.
\end{proof}

\begin{remark}
	Noting that with the regularity obtained by $\eqref{eq3.32}$-$\eqref{eq3.33}$, it follows that the problem $\eqref{eq1.01} - \eqref{eq1.04}$ has a unique strong solution $u$ satisfying
\begin{equation*}
	\left\{
	\begin{array}{l}
	u\in C^{0}\left( 0,T;H^{1}(\Omega)\right) \cap C^{1}\left( 0,T;L^{2}(\Omega)\right) \cap L^{\infty }\left( 0,T;H^{2}(\Omega)\right) \text{, }\\
	u_{t}\in L^{\infty }\left(0,T;H^{1}(\Omega)\right), u_{tt}\in L^{\infty }\left( 0,T;L^{2}(\Omega)\right) \text{, \ }\\
	u(0,\cdot )\text{,\ }u(1,\cdot )\in H^{2}\left( 0,T\right).
	\end{array}
	\right.
\end{equation*}
\end{remark}

\section{The regularity of solutions}
In this section, we study the regularity of solution of problem $\eqref{eq1.01} - \eqref{eq1.04}$. For this purpose, we also assume that the constants $h_{0}$, $h_{1}$, $\lambda _{0}$, $\lambda _{1}$, $K$, $\lambda $, $\widetilde{h}_{0}$, $\widetilde{h}_{1}$, $\widetilde{\lambda }_{0}$, $\widetilde{\lambda }_{1}$ satisfy the conditions as in section 3. Furthermore, we will impose the following stronger assumptions, with $r\in \mathbb{N}$.
\begin{description}
\item[(A1)] $\widetilde{u}_{0}\in H^{r+2}(\Omega) \, \text{ and }\, \widetilde{u}_{1}\in H^{r+1}(\Omega)$.
\item[(A2)] The function $f(x,t)$ satisfies
\begin{equation*}
	\frac{\partial^{r}f}{\partial \,x\,^{j}\partial\,t\,^{r-j}}\in L^{\infty }(0,T;L^{2}(\Omega))\text{, }0\leq j\leq r\text{,}
\end{equation*}
\begin{equation*}
	\frac{\partial \,^{\nu }f}{\partial \,t\,^{\nu }}\in L^{2}(0,T;L^{2}(\Omega))\text{, \ }0\leq \,\nu \leq r+1\text{,}
\end{equation*}
and
\begin{equation*}
	\frac{\partial \,^{\mu }f}{\partial \,t\,^{\mu }}(\cdot,0)\in H^{1}(\Omega)\text{, }0\leq \,\mu \leq r-1\text{.}
\end{equation*}
\item[(A3)] $g_{0}\text{, }g_{1}\in H^{r+1}\left( 0,T\right) \text{, }r\geq 1\text{.}$
\end{description}

Formally differentiating problem$\eqref{eq1.01} - \eqref{eq1.04}$ with respect to time up to order $r$ and letting $u^{[r]}=\frac{\partial \,^{r}u}{\partial\,t\,^{r}}$ we are led to consider the solution $u^{[r]}$ of problem $(Q^{[r]})$:
\begin{equation}\label{proQr}
	\left\{
	\begin{array}{l}
	Lu^{[r]}=f^{[r]}(x,t)\text{, \ }(x,t)\in Q_{T}\text{,\ }\\
	B_{0}u^{[r]}=g_{0}^{[r]}(t)\text{, \ }B_{1}u^{[r]}=g_{1}^{[r]}(t)\text{,}\\
	u^{[r]}(x,0)=u_{0}^{[r]}(x)\text{, \ }u_{t}^{[r]}(x,0)=u_{1}^{[r]}(x)\text{,}
	\end{array}
	\right. \tag{$Q^{[r]}$}
\end{equation}
where
\begin{equation*}
	Lu^{[r]} =u^{[r]}_{tt} - u^{[r]}_{xx} + Ku^{[r]}+\lambda u^{[r]}_{t},
\end{equation*}
\begin{equation*}
	B_{0}u^{[r]}=u^{[r]}_{x}(0,t) - h_{0}u^{[r]}(0,t) - \lambda_{0}u^{[r]}_{t}(0,t) - \widetilde{h}_{1}u^{[r]}(1,t)-\widetilde{\lambda}_{1}u^{[r]}_{t}(1,t),
\end{equation*}
\begin{equation*}
	B_{1}u^{[r]}=-u^{[r]}_{x}(1,t) - h_{1}u^{[r]}(1,t) - \lambda_{1}u^{[r]}_{t}(1,t) - \widetilde{h}_{0}u^{[r]}(0,t)-\widetilde{\lambda}_{0}u^{[r]}_{t}(0,t),
\end{equation*}
the functions $\,u_{0}^{[r]}$ and $u_{1}^{[r]}$ are defined by the recurrence formulas
\begin{equation*}
	u_{0}^{[0]}=\widetilde{u}_{0}\text{, \ }\,u_{0}^{[r]}=u_{1}^{[r-1]}\text{,\ }r\geq 1\text{, }
\end{equation*}
\begin{equation*}
	u_{1}^{[0]}=\widetilde{u}_{1}\text{, \ }u_{1}^{[r]}=u_{0xx}^{[r-1]}\,-Ku_{0}^{[r-1]}-\lambda u_{1}^{[r-1]}+\frac{\partial \,^{r-1}f}{\partial \,t\,^{r-1}}(x,0)\text{, \ }r\geq 1\text{,}
\end{equation*}
and
\begin{equation*}
	f^{[r]}=\frac{\partial \,^{r}f}{\partial \,t\,^{r}}, \quad g_{i}^{[0]}=g_{i}\text{, \ }g_{i}^{[r]}=\frac{d\,^{r}g_{i}}{d\,t\,^{r}}\text{, }r\geq 1\text{, }i=0,1\text{.\ }
\end{equation*}
From the assumptions (A1)-(A3) we deduce that $u_{0}^{[r]}$, $u_{1}^{[r]}$, $f^{[r]}$, $g_{0}^{[r]}$ and $g_{1}^{[r]}$ satisfy the conditions of Theorem 3.2. Hence, the problem $(Q^{[r]})$ has a unique weak solution $u^{[r]}$ such that
\begin{equation}\label{eq4.01}
	\left\{
\begin{array}{l}
u^{[r]}\in L^{\infty }\left( 0,T;H^{2}(\Omega)\right) \cap C^{0}\left(0,T;H^{1}(\Omega)\right) \cap C^{1}\left( 0,T;L^{2}(\Omega)\right), \\
u_{t}^{[r]}\in L^{\infty }\left( 0,T;H^{1}(\Omega)\right), \, u_{tt}^{[r]}\in L^{\infty }\left( 0,T;L^{2}(\Omega)\right), \\
u^{[r]}(0,\cdot )\text{, \ }u^{[r]}(1,\cdot )\in H^{2}\left( 0,T\right) \text{.}%
\end{array}%
\right. \tag{4.1}
\end{equation}
Moreover, from the uniqueness of weak solution we have $\,u^{[r]}=\frac{\partial \,^{r}u}{\partial \,t\,^{r}}$. Hence we deduce from $\eqref{eq4.01}$ that the solution $u$ of problem $\eqref{eq1.01}-\eqref{eq1.04}$ satisfy
\begin{gather}
u\in C^{r-1}\left( 0,T;H^{2}(\Omega)\right) \cap C^{r}\left( 0,T;H^{1}(\Omega)\right) \cap C^{r+1}\left( 0,T;L^{2}(\Omega)\right) \text{,} \notag \\
\frac{\partial \,^{r}u}{\partial \,t\,^{r}}\in L^{\infty }\left(0,T;H^{2}(\Omega)\right) \cap C^{0}\left( 0,T;H^{1}(\Omega)\right) \cap C^{1}\left(0,T;L^{2}(\Omega)\right) \text{, } \notag \\
\frac{\partial \,^{r+1}u}{\partial \,t\,^{r+1}}\in L^{\infty }\left(0,T;H^{1}(\Omega)\right) \text{, } \label{eq4.02} \tag{4.2}\\
\frac{\partial \,^{r+2}u}{\partial \,t\,^{r+2}}\in L^{\infty }\left(0,T;L^{2}(\Omega)\right) \text{, } \notag\\
u(0,\cdot )\text{, }\,u(1,\cdot )\in H^{r+2}\left( 0,T\right) \text{.} \notag %
\end{gather}
\indent Next we shall prove by induction on $r$ that
\begin{equation}\label{eq4.03}
	\frac{\partial \,^{r+2-j}u}{\partial \,t\,^{r+2-j}}\in L^{\infty }(0,T;H^{j}(\Omega))\text{,\ \ }0\leq j\leq r+2\text{.} \tag{4.3}
\end{equation}
With $r=1$, it follows from $\eqref{eq4.02}$ that
\begin{equation}\label{eq4.04}
	u'\in L^{\infty }(0,T;H^{2}(\Omega)), \, u'' \in L^{\infty }(0,T;H^{1}(\Omega)), \, u'''\in L^{\infty }(0,T;L^{2}(\Omega)). \tag{4.4}
\end{equation}
On the other hand, from $\eqref{eq1.01}$, $\eqref{eq4.04}$ and the assumption (A2) we deduce that
\begin{equation*}
	u_{xxx} = u''_{x} + Ku_{x} +\lambda u'_{x} - f_{x} \in L^{\infty }(0,T;L^{2}(\Omega)).
\end{equation*}
Thus, $u \in L^{\infty }(0,T;H^{3}(\Omega))$ and $\eqref{eq4.03}$ hold for $r =1$. Suppose by induction that $\eqref{eq4.03}$ holds for $r-1$, i.e.,
\begin{equation}\label{eq4.05}
	\frac{\partial \,^{r+1-j}u}{\partial \,t\,^{r+1-j}}\in L^{\infty }(0,T;H^{j}(\Omega))\text{,\ \ }0\leq j\leq r+1\text{.} \tag{4.5}
\end{equation}
We shall prove that $\eqref{eq4.03}$ holds. It follows from $\eqref{eq4.02}$ that
\begin{equation}\label{eq4.06}
	\frac{\partial \,^{r+2-j}u}{\partial \,t\,^{r+2-j}}\in L^{\infty }\left(0,T;H^{j}(\Omega)\right) \text{,} \, \text{for} \,\, j =0, 1, 2. \tag{4.6}
\end{equation}
Let $j\in \left\{3, 4, ..., r+2\right\}$ and put $\theta  = r+2 -j$. We have from $\eqref{eq1.01}$
\begin{equation}\label{eq4.07}
	\frac{\partial^{r+2}u}{\partial x^j \partial t^{\theta}} = \frac{\partial^{r+2}u}{\partial x^{j-2} \partial t^{\theta +2}} + K\frac{\partial^{r}u}{\partial x^{j-2} \partial t^{\theta}}+ \lambda \frac{\partial^{r+1}u}{\partial x^{j-2} \partial t^{\theta+1}} - \frac{\partial^{r}f}{\partial x^{j-2} \partial t^{\theta}}. \tag{4.7}
\end{equation}
On the other hand, it follows from $\eqref{eq4.05}$ and the assumption (A2), that
\begin{equation}\label{eq4.08}
	\frac{\partial^{\theta}u}{\partial t^{\theta}}\in L^{\infty }(0,T;H^{j-1}(\Omega)), \, \frac{\partial^{\theta+1}u}{\partial t^{\theta+1}}, \, \, \frac{\partial^{\theta}f}{\partial t^{\theta}} \in L^{\infty }(0,T;H^{j-2}(\Omega)). \tag{4.8}
\end{equation}
Combining $\eqref{eq4.06}$, $\eqref{eq4.07}$ and $\eqref{eq4.08}$, by induction arguments on $j$, we conclude that $\eqref{eq4.03}$ holds.

\vspace{0.2cm}
Hence we have the following theorem
\begin{theorem}
	Let (A1)-(A3) hold. Then the unique solution $u(x,t)$ of problem $\eqref{eq1.01}$-$\eqref{eq1.04}$ satisfies $\eqref{eq4.02}$ and $\eqref{eq4.03}.$ Furthermore
\begin{equation}\label{eq4.09}
	u\in H^{r+2}\left(Q_{T}\right) \cap \left( \bigcap_{j=0}^{r+1}C^{r+1-j}\left(0,T; H^{j}(\Omega)\right) \right). \tag{4.9}
\end{equation}
\end{theorem}

\section{Exponential decay of solutions}
	In this section we assume that $K>0$ and $\lambda>0$. Let $u(x,t)$ be a strong solution of problem $\eqref{eq1.01}$-$\eqref{eq1.04}$. In order to obtain the decay result, we use the functional
\begin{equation}\label{eq5.01}
	\Gamma(t) = E(t) + \delta \psi(t), \tag{5.1}
\end{equation}
where $\delta$ is a positive constant and
\begin{equation}\label{eq5.02}
	E\left( t\right) =\frac{1}{2}\left\Vert u^{\prime }\left( t\right)\right\Vert ^{2}+\frac{1}{2}\left\Vert u\left( t\right) \right\Vert _{a}^{2} + \frac{K}{2} \left\|u(t)\right\|^2, \tag{5.2}
\end{equation}
\begin{equation}\label{eq5.03}
	\psi \left( t\right) =\left\langle u\left( t\right) ,u^{\prime }\left(t\right) \right\rangle +\frac{\lambda }{2}\left\Vert u\left( t\right) \right\Vert ^{2}+\frac{\lambda _{0}}{2} u^{2}\left( 0,t\right) +\frac{\lambda _{1}}{2}u^{2}\left( 1,t\right). \tag{5.3}
\end{equation}
\begin{lemma}
	There exist the constants $\beta_{1}$, $\beta_{2}$ such that
\begin{equation}\label{eq5.04}
	\beta _{1}E\left( t\right) \leq \Gamma \left( t\right) \leq \beta_{2}E\left( t\right), \tag{5.4}
\end{equation}
when $\delta < \frac{C_0}{2}$.
\end{lemma}
\begin{proof}
	By using Lemma 2.1, it's easy to obtain the following estimate
\begin{equation*}
	\Gamma (t)\leq \frac{1+\delta }{2}\left\Vert u^{\prime }(t)\right\Vert ^{2}+%
\frac{1}{2}\left[ 1+\frac{2\delta }{C_{0}}\left( 1+\lambda +\lambda
_{0}+\lambda _{1}\right) \right] \left\Vert u(t)\right\Vert _{a}^{2}+\frac{K%
}{2}\left\Vert u(t)\right\Vert ^{2}\text{,}%
\end{equation*}
which implies that
\begin{equation*}
	\Gamma \left( t\right) \leq \beta_{2}E\left( t\right),
\end{equation*}
where
\begin{equation*}
	\beta _{2}=1+\frac{2\delta }{C_{0}}\left( 1+\lambda +\lambda _{0}+\lambda
_{1}\right) \text{.}%
\end{equation*}
Similar, we have
\begin{equation*}
	\Gamma (t)\geq \frac{1-\delta }{2}\left\Vert u^{\prime }(t)\right\Vert
^{2}+\left( \frac{1}{2}-\frac{\delta }{C_{0}}\right) \left\Vert
u(t)\right\Vert _{a}^{2}+\frac{K}{2}\left\Vert u(t)\right\Vert ^{2}\text{.}%
\end{equation*}
Thus, if $\delta <\frac{C_{0}}{2}$ then $\Gamma(t) \geq \beta_{1}E(t)$, where $\beta _{1}=1-\frac{2\delta }{C_{0}}>0$. Lemma 5.1 is proved.
\end{proof}

\begin{lemma}
	The functional $E(t)$ defined by \eqref{eq5.02}, satisfies
\begin{gather}\label{eq5.05}
	E^{\prime }(t)\leq (\frac{\varepsilon _{1}}{2}-\lambda )\left\Vert u^{\prime
}(t)\right\Vert ^{2}+\left( \varepsilon _{1}-\frac{\mu _{\min }}{2}\right) %
\left[ \left\vert u^{\prime }(0,t)\right\vert ^{2}+\left\vert u^{\prime
}(1,t)\right\vert ^{2}\right] \tag{5.5} \\
\text{ \ \ \ \ \ }+\frac{1}{\varepsilon _{1}C_{0}}\left( \widetilde{h}%
_{0}^{2}+\widetilde{h}_{1}^{2}\right) \left\Vert u(t)\right\Vert _{a}^{2}+%
\frac{1}{2\varepsilon _{1}}\left[ g_{0}^{2}(t)+g_{1}^{2}(t)+\left\Vert
f(t)\right\Vert ^{2}\right] \text{,}%
\notag
\end{gather}
for all $\varepsilon_{1}>0 $.
\end{lemma}
\begin{proof}
	Multiplying $\eqref{eq1.01}$ by $u'(x,t)$ and integrating over $[0,1]$, we get
\begin{alignat}{2}\label{eq5.06}
E'\left( t\right)  =&-\lambda \left\Vert u^{\prime }\left(t\right) \right\Vert ^{2}-\left\{ \lambda _{0}\left\vert u^{\prime }\left(0,t\right) \right\vert ^{2}+\lambda _{1}\left\vert u^{\prime }\left(1,t\right) \right\vert ^{2}\right.   \notag \\
&\left. +\left( \widetilde{\lambda }_{0}+\widetilde{\lambda }_{1}\right)u^{\prime }\left( 0,t\right) u^{\prime }\left( 1,t\right)^{{}}\right\} -\widetilde{h}_{0}u\left( 0,t\right) u^{\prime }\left( 1,t\right)  \tag{5.6} \\
&-\widetilde{h}_{1}u\left( 1,t\right) u^{\prime }\left( 0,t\right)-g_{0}\left( t\right) u^{\prime }\left( 0,t\right) -g_{1}\left( t\right)u^{\prime }\left( 1,t\right)   \notag
&+\left\langle f\left( t\right) ,u^{\prime }\left( t\right) \right\rangle . \notag
\end{alignat}
By Lemma 1.3 we have
\begin{gather}
\lambda _{0}\left\vert u^{\prime }\left( 0,t\right) \right\vert ^{2}+\lambda_{1}\left\vert u^{\prime }\left( 1,t\right) \right\vert ^{2}+\left( \widetilde{\lambda }_{0}+\widetilde{\lambda }_{1}\right) u^{\prime }\left(0,t\right) u^{\prime }\left( 1,t\right)   \label{eq5.07} \tag{5.7}\\
\geq \frac{\mu _{\min }}{2}\left[ \left\vert u^{\prime }\left( 0,t\right)\right\vert ^{2}+\left\vert u^{\prime }\left( 1,t\right) \right\vert ^{2}\right].   \notag
\end{gather}
It follows from $\eqref{eq5.06}$ and $\eqref{eq5.07}$ that
\begin{alignat}{2}\label{eq5.08}
E'\left( t\right)  \leq &-\lambda \left\Vert u^{\prime }\left(t\right) \right\Vert ^{2}-\frac{\mu _{\min }}{2}\left[ \left\vert u^{\prime}\left( 0,t\right) \right\vert ^{2}+\left\vert u^{\prime }\left( 1,t\right)\right\vert ^{2}\right]   \notag \\
&-\widetilde{h}_{0}u\left( 0,t\right) u^{\prime }\left( 1,t\right) -\widetilde{h}_{1}u\left( 1,t\right) u^{\prime }\left( 0,t\right)   \tag{5.8}\\
&-g_{0}\left( t\right) u^{\prime }\left( 0,t\right) -g_{1}\left( t\right)u^{\prime }\left( 1,t\right) +\left\langle f\left( t\right) ,u^{\prime}\left( t\right) \right\rangle.   \notag
\end{alignat}
On the other hand, for $\varepsilon_{1}>0$,
\begin{alignat}{2}\label{eq5.9}
-\widetilde{h}_{0}u(0,t)u^{\prime }(1,t)&\leq \frac{\varepsilon _{1}}{2}%
\left\vert u^{\prime }(1,t)\right\vert ^{2}+\frac{1}{2\varepsilon _{1}}%
\widetilde{h}_{0}^{2}u^{2}(0,t)   \tag{5.9} \\
&\leq \frac{\varepsilon _{1}}{2}\left\vert u^{\prime
}(1,t)\right\vert ^{2}+\frac{1}{\varepsilon _{1}C_{0}}\widetilde{h}%
_{0}^{2}\left\Vert u(t)\right\Vert _{a}^{2}\text{,}
 \notag
\end{alignat}
\begin{equation}\label{eq5.10}
-\widetilde{h}_{1}u(1,t)u^{\prime }(0,t)\leq \frac{\varepsilon _{1}}{2}%
\left\vert u^{\prime }(0,t)\right\vert ^{2}+\frac{1}{\varepsilon _{1}C_{0}}%
\widetilde{h}_{1}^{2}\left\Vert u(t)\right\Vert _{a}^{2}\text{,} \tag{5.10}
\end{equation}
\begin{alignat}{2}\label{eq5.11}
	-g_{0}(t)u^{\prime }(0,t)-g_{1}(t)u^{\prime }(1,t)  & \leq   \frac{\varepsilon
_{1}}{2}\left[ \left\vert u^{\prime }(0,t)\right\vert ^{2}+\left\vert
u^{\prime }(1,t)\right\vert ^{2}\right] \tag{5.11} \\
	&+ \frac{1}{2\varepsilon _{1}}\left[g_{0}^{2}(t)+g_{1}^{2}(t)\right], \notag
\end{alignat}
\begin{equation}\label{eq5.12}
\langle f(t),u^{\prime }(t)\rangle \leq \left\Vert f(t)\right\Vert
\left\Vert u^{\prime }(t)\right\Vert \leq \frac{\varepsilon _{1}}{2}%
\left\Vert u^{\prime }(t)\right\Vert ^{2}+\frac{1}{2\varepsilon _{1}}%
\left\Vert f(t)\right\Vert ^{2}.  \tag{5.12}
\end{equation}
Combining $\eqref{eq5.08}$ - $\eqref{eq5.12}$, it is easy to see that $\eqref{eq5.05}$ holds. The proof is complete.
\end{proof}
\begin{lemma}
	The functional $\psi(t)$ defined by $\eqref{eq5.03}$ satisfies
\begin{gather}
\psi ^{\prime }(t)\leq \left\Vert u^{\prime }(t)\right\Vert ^{2}+\left(
\frac{2}{C_{0}}\left\vert \widetilde{h}_{0}+\widetilde{h}_{1}\right\vert +%
\frac{5\varepsilon _{2}}{C_{0}}-1\right) \left\Vert u(t)\right\Vert _{a}^{2} \notag \\
+\frac{1}{2\varepsilon _{2}}\left( \widetilde{\lambda
}_{0}^{2}+\widetilde{\lambda }_{1}^{2}\right) \left[ \left\vert u^{\prime
}(0,t)\right\vert ^{2}+\left\vert u^{\prime }(1,t)\right\vert ^{2}\right]   \label{eq5.13} \tag{5.13} \\
+\frac{1}{2\varepsilon _{2}}\left[ \left\Vert
f(t)\right\Vert ^{2}+g_{0}^{2}(t)+g_{1}^{2}(t)\right] \text{,} \notag
\end{gather}
for all $\varepsilon_{2}>0$.
\end{lemma}
\begin{proof}
	Multiplying the equation $\eqref{eq1.01}$ by $u(x,t)$ and integrating over $[0,1]$, we have
\begin{gather}
\psi ^{\prime }(t)=\left\Vert u^{\prime }(t)\right\Vert ^{2}-\left\Vert
u(t)\right\Vert _{a}^{2}-K\left\Vert u(t)\right\Vert ^{2}-\left( \widetilde{h%
}_{0}+\widetilde{h}_{1}\right) u(0,t)u(1,t)  \notag \\
-\widetilde{\lambda }_{0}u^{\prime }(0,t)u(1,t)-%
\widetilde{\lambda }_{1}u(0,t)u^{\prime }(1,t)  \label{eq5.14} \tag{5.14}\\
 -g_{0}(t)u(0,t)-g_{1}(t)u(1,t)+\langle
f(t),u(t)\rangle \text{.}   \notag
\end{gather}
By some estimations as in proof of Lemma 5.2, we deduce the conclusion of Lemma.
\end{proof}
\begin{theorem} Assume that
\begin{equation}\label{eq5.15}
	\sigma(t) \leq \sigma_{1}\exp(-\sigma_{2}t), \, \text{for all} \, \, t \geq 0, \tag{5.15}
\end{equation}
where $\sigma_{1}$, $\sigma_{2}$ are two positive constants and
\begin{equation*}
	\sigma(t) =  \left\|f(t)\right\|^{2}+ g_{0}^{2}(t) +g_{1}^{2}(t).
\end{equation*}
Then, there exist positive constants $\gamma_{1}$, $\gamma_{2}$ such that
\begin{equation}\label{eq5.16}
	E(t) \leq \gamma_{1}\exp(-\gamma_{2}t), \, \text{for all} \,\, t \geq 0, \tag{5.16}
\end{equation}
for any strong solution of the problem $\eqref{eq1.01}$-$\eqref{eq1.04},$ where $\widetilde{h}_{0}$ and $\widetilde{h}_{1}$ are chosen small enough.
\end{theorem}
\begin{proof}
	It follows from $\eqref{eq5.01}$, $\eqref{eq5.05}$ and $\eqref{eq5.13}$, that
\begin{alignat}{2}\label{eq5.17}
\Gamma ^{\prime }(t)\leq & \left( \delta +\frac{\varepsilon _{1}}{C_{0}}%
-\lambda \right) \left\Vert u^{\prime }(t)\right\Vert ^{2} \notag\\
&+\left[ \frac{1}{\varepsilon _{1}C_{0}}\left( \widetilde{h}%
_{0}^{2}+\widetilde{h}_{1}^{2}\right) +\delta \left( \frac{2}{C_{0}}%
\left\vert \widetilde{h}_{0}+\widetilde{h}_{1}\right\vert +\frac{%
5\varepsilon _{2}}{C_{0}}-1\right) \right] \left\Vert u(t)\right\Vert
_{a}^{2} \notag \\
&+\left[ \varepsilon _{1}-\frac{\mu _{\min }}{2}+\frac{\delta
}{2\varepsilon _{2}}\left( \widetilde{\lambda }_{0}^{2}+\widetilde{\lambda }%
_{1}^{2}\right) \right] \left[ \left\vert u^{\prime }(0,t)\right\vert
^{2}+\left\vert u^{\prime }(1,t)\right\vert ^{2}\right] \tag{5.17}\\
&+\frac{1}{2}\left( \frac{1}{\varepsilon _{1}}+\frac{\delta }{%
\varepsilon _{2}}\right) \left[ \left\Vert f(t)\right\Vert
^{2}+g_{0}^{2}(t)+g_{1}^{2}(t)\right] \text{,} \notag
\end{alignat}
for all $\varepsilon_{1}$, $\varepsilon_{2} > 0$.

\vspace{0.3cm}
Let $\varepsilon _{1}<\min \{C_{0}\lambda ,\frac{\mu _{\min }}{2}\}$, $%
\varepsilon _{2}<\frac{1}{5} C_{0}$ and
\begin{equation*}
\delta <\min \Big \{\frac{C_{0}}{2}\text{, }\lambda -\frac{\varepsilon _{1}}{C_{0}%
}\text{, }\frac{2\varepsilon _{2}}{\widetilde{\lambda }_{0}^{2}+\widetilde{%
\lambda }_{1}^{2}}\left( \frac{\mu _{\min }}{2}-\varepsilon _{1}\right) \Big \}.%
\end{equation*}
Then, by choosing $\widetilde{h}_{0}$, $\widetilde{h}_{1}$ satisfy
\begin{equation*}
\frac{1}{\varepsilon _{1}C_{0}}\left( \widetilde{h}_{0}^{2}+\widetilde{h}%
_{1}^{2}\right) +\frac{2\delta }{C_{0}}\left\vert \widetilde{h}_{0}+%
\widetilde{h}_{1}\right\vert <\delta \left( 1-\frac{5\varepsilon _{2}}{C_{0}}%
\right) \text{,}%,
\end{equation*}
 we deduce from $\eqref{eq5.04}$ and $\eqref{eq5.17}$ that there exists a constant $\gamma < \sigma_{2}$ such that
\begin{equation}\label{eq5.18}
\Gamma ^{\prime }(t)\leq -\delta \Gamma (t)+\frac{1}{2}\left( \frac{1}{%
\varepsilon _{1}}+\frac{\delta }{\varepsilon _{2}}\right) \sigma (t)	, \text{\ for all} \,\, t\geq 0. \tag{5.18}
\end{equation}
Combining $\eqref{eq5.04}$, $\eqref{eq5.15}$ and $\eqref{eq5.18}$, we get $\eqref{eq5.16}$. Theorem 5.4 is completely proved.
\end{proof}

We can extend the above theorem to weak solutions by using density arguments.

\section{\textbf{Numerical results}}

$\qquad $Consider the following problem:{\small
\begin{equation}
u_{tt}-u_{xx}+Ku+\lambda u_{t}=f(x,t)\hspace{0.2cm}\text{in}\hspace{0.2cm}%
(0,1)\times (0,\infty ),  \tag{6.1}  \label{s1}
\end{equation}%
}\newline
$0<x<1,$ $0<t<T,$ with boundary conditions{\small
\begin{equation}
\begin{tabular}{l}
$u_{x}(0,t)=h_{0}u(0,t)+\lambda _{0}u_{t}(0,t)+\widetilde{h}_{1}u(1,t)+%
\widetilde{\lambda }_{1}u_{t}(1,t)+g_{0}(t),$\textit{\medskip } \\
$-u_{x}(1,t)=h_{1}u(1,t)+\lambda _{1}u_{t}(1,t)+\widetilde{h}_{0}u(0,t)+%
\widetilde{\lambda }_{0}u_{t}(0,t)+g_{1}(t)\text{,}$%
\end{tabular}
\tag{6.2}  \label{s2}
\end{equation}%
}$\newline
$and initial conditions{\small
\begin{equation}
u(x,0)=\widetilde{u}_{0}(x),\text{ \ }u_{t}(x,0)=\widetilde{u}_{1}(x),
\tag{6.3}  \label{s3}
\end{equation}%
}$\newline
$where $K=\lambda =\lambda _{0}=\lambda _{1}=1,$ $\widetilde{\lambda }_{0}=%
\widetilde{\lambda }_{1}=\frac{-1}{2},$ $h_{0}=h_{1}=1,$ $\widetilde{h}_{0}=%
\frac{1}{2},$ $\widetilde{h}_{1}=\frac{-3}{2}$ are constants and the
functions $\widetilde{u}_{0},$ $\widetilde{u}_{1},$ $g_{0},$\ $g_{1}$ and $f$
are defined by%
\begin{equation}
u_{0}(x)=1+(x^{2}-x)^{2},\text{\ }\widetilde{u}_{1}(x)=-1-(x^{2}-x)^{2},
\tag{6.4}  \label{s4}
\end{equation}%
\begin{equation}
{\small g_{0}(t)=e}^{-t},\text{ }{\small g_{1}(t)=-{\small e}^{-t}},
\tag{6.5}  \label{s5}
\end{equation}%
\begin{equation}
f(x,t)=\left( x^{4}-2x^{3}-11x^{2}-12x-1\right) e^{-t}.  \tag{6.6}
\label{s6}
\end{equation}%
$\qquad $The exact solution of the problem (\ref{s1}) -- (\ref{s3}) with $%
\widetilde{u}_{0},$ $\widetilde{u}_{1},$ $g_{0},$\ $g_{1}$ and $f$ defined
in (\ref{s4}) -- (\ref{s6}) respectively, is the function $U_{ex}$ given by%
{\small
\begin{equation}
U_{ex}(x,t)=\left( x^{4}-2x^{3}+x^{2}+1\right) e^{-t}.  \tag{6.7}  \label{s7}
\end{equation}%
}

$\qquad $To solve problem (\ref{s1}) -- (\ref{s3}) numerically, we consider
the differential system for the unknowns $U_{j}(t)\equiv u(x_{j},t),$ $%
V_{j}(t)=\frac{dU_{j}}{dt}(t),$ with $x_{j}=j\Delta x,$ $\Delta x=\frac{1}{N}%
,$ $j=0,1,...,N:$

{\small
\begin{equation}
\left\{
\begin{tabular}{l}
$\frac{dU_{j}}{dt}(t)=V_{j}(t),\text{\ }j=\overline{0,N},$\textit{\medskip }
\\
$\frac{dV_{0}}{dt}(t)=-\left( \frac{1+h_{0}\Delta x}{\left( \Delta x\right)
^{2}}+K\right) U_{0}(t)+\frac{1}{\left( \Delta x\right) ^{2}}U_{1}(t)-\frac{%
\widetilde{h}_{1}}{\Delta x}U_{N}(t)$\textit{\medskip } \\
$\ \ \ \ \ \ \ \ \ \ \ \ \ \ \ -(\frac{\lambda _{0}}{\Delta x}+\lambda
)V_{0}(t)-\frac{\widetilde{\lambda }_{1}}{\Delta x}V_{N}(t)-\frac{1}{\Delta x%
}g_{0}(t)+f_{0}(t),$ \\
$\frac{dV_{j}}{dt}(t)=\frac{U_{j-1}(t)-2U_{j}(t)+U_{j+1}(t)}{\left( \Delta
x\right) ^{2}}-KU_{j}(t)-\lambda V_{j}(t)+f(x_{j},t),$\ $j=\overline{1,N-1},$%
\textit{\medskip } \\
$\frac{dV_{N}}{dt}(t)=-\frac{\widetilde{h}_{0}}{\Delta x}U_{0}(t)+\frac{1}{%
\left( \Delta x\right) ^{2}}U_{N-1}(t)-\left( \frac{1+h_{1}\Delta x}{\left(
\Delta x\right) ^{2}}+K\right) U_{N}(t)$\textit{\medskip } \\
$\ \ \ \ \ \ \ \ \ \ \ \ \ \ \ -\frac{\widetilde{\lambda }_{0}}{\Delta x}%
V_{0}(t)-(\frac{\lambda _{1}}{\Delta x}+\lambda )V_{N}(t)-\frac{1}{\Delta x}%
g_{1}(t)+f_{N}(t),$\textit{\medskip } \\
$U_{j}(0)=\widetilde{u}_{0}(x_{j}),\text{ }V_{j}(0)=\widetilde{u}_{1}(x_{j}),%
\text{\ }j=\overline{0,N}.$\textit{\medskip }%
\end{tabular}%
\right.   \tag{6.8}  \label{s8}
\end{equation}%
}\newline
Then system (\ref{s8}) is equivalent to:$\newline
${\small
\begin{equation}
\frac{d}{dt}\left[
\begin{tabular}{l}
$U_{0}$ \\
$U_{1}$ \\
$\vdots $ \\
$\vdots $ \\
$U_{N}$ \\
$V_{0}$ \\
$V_{1}$ \\
$\vdots $ \\
$\vdots $ \\
$V_{N}$%
\end{tabular}%
\right] =\left[
\begin{tabular}{|lllll|lllll|}
\hline
0 &  &  &  &  & $1$ &  &  &  &  \\
0 & 0 &  &  &  &  & $1$ &  &  &  \\
& $\ddots $ & $\ddots $ &  &  &  &  & $\ddots $ &  &  \\
&  & $\ddots $ & $\ddots $ &  &  &  &  & $\ddots $ &  \\
&  &  & 0 & 0 &  &  &  &  & $1$ \\ \hline
\multicolumn{1}{|r}{$\widetilde{\gamma }_{0}$} & \multicolumn{1}{r}{$\alpha $%
} &  &  & $\widetilde{\gamma }_{N}$ & $\widetilde{\delta }_{0}$ &  &  &  & $%
\widetilde{\delta }_{N}$ \\
\multicolumn{1}{|r}{$\alpha $} & \multicolumn{1}{r}{$\gamma $} & $\alpha $ &
&  &  & $-\lambda $ &  &  &  \\
& \multicolumn{1}{c}{$\ddots $} & \multicolumn{1}{c}{$\ddots $} &
\multicolumn{1}{c}{$\ddots $} &  &  &  & $\ddots $ &  &  \\
&  & \multicolumn{1}{r}{$\alpha $} & \multicolumn{1}{r}{$\gamma $} & $\alpha
$ &  &  &  & $-\lambda $ &  \\
$\widehat{\gamma }_{0}$ &  &  & \multicolumn{1}{r}{$\alpha $} &
\multicolumn{1}{r|}{$\widehat{\gamma }_{N}$} & $\widehat{\delta }_{0}$ &  &
&  & $\widehat{\delta }_{N}$ \\ \hline
\end{tabular}%
\right] \left[
\begin{tabular}{l}
$U_{0}$ \\
$U_{1}$ \\
$\vdots $ \\
$\vdots $ \\
$U_{N}$ \\
$V_{0}$ \\
$V_{1}$ \\
$\vdots $ \\
$\vdots $ \\
$V_{N}$%
\end{tabular}%
\right] +\left[
\begin{tabular}{l}
$0$ \\
$0$ \\
$\vdots $ \\
$\vdots $ \\
$0$ \\
$F_{0}$ \\
$F_{1}$ \\
$\vdots $ \\
$\vdots $ \\
$F_{N}$%
\end{tabular}%
\right]   \tag{6.9}  \label{s9}
\end{equation}%
}\newline
where\newline
{\small
\begin{equation}
\left\{
\begin{tabular}{l}
$\alpha =\frac{1}{\left( \Delta x\right) ^{2}},\ \gamma =-K-\frac{2}{\left(
\Delta x\right) ^{2}}=-K-2\alpha ,$\textit{\medskip } \\
$\widetilde{\gamma }_{0}=-\left( K+\frac{1+h_{0}\Delta x}{\left( \Delta
x\right) ^{2}}\right) ,$ $\widetilde{\gamma }_{N}=-\frac{\widetilde{h}_{1}}{%
\Delta x},$ $\widetilde{\delta }_{0}=-(\frac{\lambda _{0}}{\Delta x}+\lambda
),$ $\widetilde{\delta }_{N}=-\frac{\widetilde{\lambda }_{1}}{\Delta x},$%
\textit{\medskip } \\
$\widehat{\gamma }_{0}=-\frac{\widetilde{h}_{0}}{\Delta x},$ $\widehat{%
\gamma }_{N}=-\left( K+\frac{1+h_{1}\Delta x}{\left( \Delta x\right) ^{2}}%
\right) ,$ $\widehat{\delta }_{0}=-\frac{\widetilde{\lambda }_{0}}{\Delta x},
$ $\widehat{\delta }_{N}=-(\frac{\lambda _{1}}{\Delta x}+\lambda ).$\textit{%
\medskip } \\
$F_{j}=F_{j}(t)=f_{j}(t)=f(x_{j},t),\ j=\overline{1,N-1},$\textit{\medskip }
\\
$F_{0}=F_{0}(t)=-\frac{1}{\Delta x}g_{0}(t)+f_{0}(t),$\textit{\medskip } \\
$F_{N}=F_{N}(t)=-\frac{1}{\Delta x}g_{1}(t)+f_{N}(t),$\textit{\medskip }%
\end{tabular}%
\right.   \tag{6.10}  \label{s10}
\end{equation}%
}\newline

 Rewritten (\ref{s9}){\small
\begin{equation}
\begin{array}{c}
\frac{d}{dt}X(t)=AX(t)+F(t)\text{,}%
\end{array}
\tag{6.11}  \label{s11}
\end{equation}%
}$\bigskip $

{\small
\begin{equation}
\left\{
\begin{tabular}{l}
$X(t)=\left(
U_{0}(t),U_{1}(t),...,U_{N}(t),V_{0}(t),V_{1}(t),...,V_{N}(t)\right) ^{T}\in
%TCIMACRO{\U{211d} }%
%BeginExpansion
\mathbb{R}
%EndExpansion
^{2N+2},$\textit{\medskip } \\
$F(t)=\left( 0,0,...,0,F_{0},F_{1},...,F_{N}\right) ^{T}\in
%TCIMACRO{\U{211d} }%
%BeginExpansion
\mathbb{R}
%EndExpansion
^{2N+2},$\textit{\medskip } \\
$A=\left[
\begin{tabular}{ll}
$O$ & $E$ \\
$\widetilde{A}$ & $\widetilde{B}$%
\end{tabular}%
\right] ,$%
\end{tabular}%
\right.  \tag{6.12}  \label{s12}
\end{equation}%
}$\bigskip $

{\small
\begin{equation}
E=\left[
\begin{tabular}{|lllll|}
\hline
$1$ &  &  &  &  \\
& $1$ &  &  &  \\
&  & $\ddots $ &  &  \\
&  &  & $\ddots $ &  \\
&  &  &  & $1$ \\ \hline
\end{tabular}%
\right] ,\ \widetilde{A}=\left[
\begin{tabular}{|rrlll|}
\hline
$\widetilde{\gamma }_{0}$ & $\alpha $ &  &  & $\widetilde{\gamma }_{N}$ \\
$\alpha $ & $\gamma $ & $\alpha $ &  &  \\
\multicolumn{1}{|l}{} & \multicolumn{1}{c}{$\ddots $} & \multicolumn{1}{c}{$%
\ddots $} & \multicolumn{1}{c}{$\ddots $} &  \\
\multicolumn{1}{|l}{} & \multicolumn{1}{l}{} & \multicolumn{1}{r}{$\alpha $}
& \multicolumn{1}{r}{$\gamma $} & $\alpha $ \\
\multicolumn{1}{|l}{$\widehat{\gamma }_{0}$} & \multicolumn{1}{l}{} &  &
\multicolumn{1}{r}{$\alpha $} & \multicolumn{1}{r|}{$\widehat{\gamma }_{N}$}
\\ \hline
\end{tabular}%
\right] ,\ \widetilde{B}=\left[
\begin{tabular}{|lllll|}
\hline
$\widetilde{\delta }_{0}$ &  &  &  & $\widetilde{\delta }_{N}$ \\
& $-\lambda $ &  &  &  \\
&  & $\ddots $ &  &  \\
&  &  & $-\lambda $ &  \\
$\widehat{\delta }_{0}$ &  &  &  & $\widehat{\delta }_{N}$ \\ \hline
\end{tabular}%
\right]  \tag{6.13}  \label{s13}
\end{equation}%
}$\qquad $To solve the linear differential system (\ref{s11}), we
use a spectral method with a time step $\Delta t=0.08 $ and a
spacial step $\Delta x=0.1 $

\begin{center}
\includegraphics[width=7cm,height=7cm,
 keepaspectratio=true]{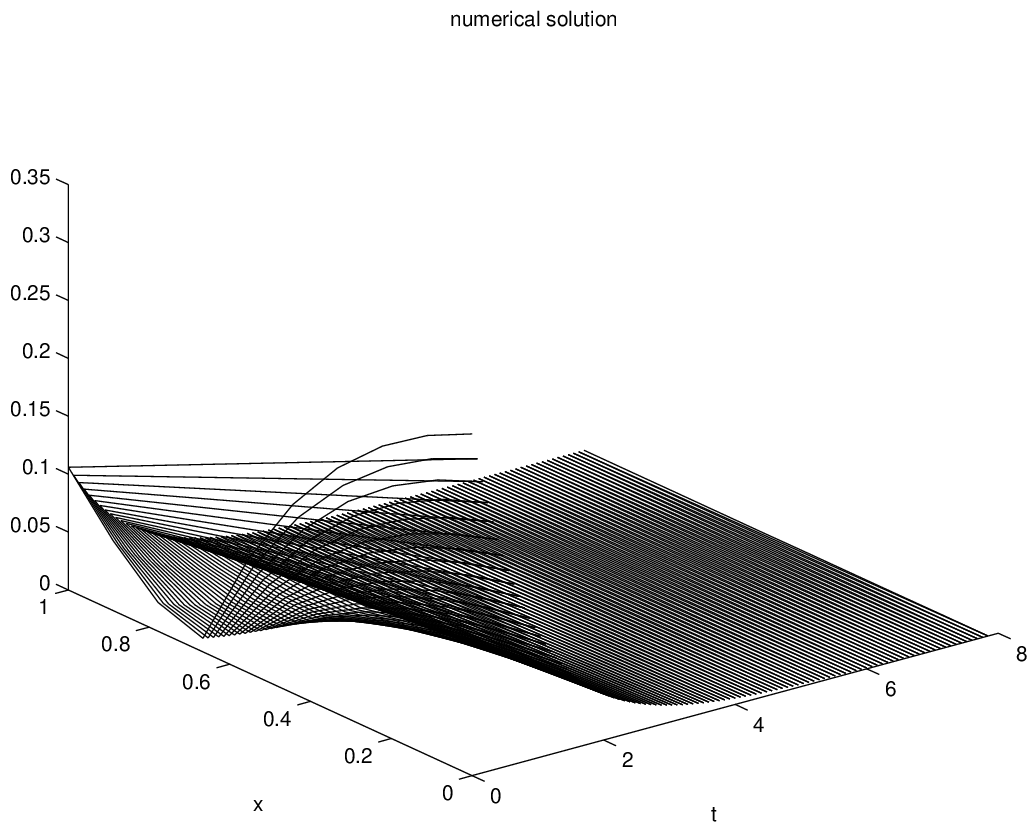}

Figure 1: Numerical solution
\end{center}

In fig.1 we have drawn the approximated solution of the problem
\ref{s1}-\ref{s3} while fig.2 represents his corresponding exact
solution \ref{s7}.

\begin{center}
\includegraphics[width=7cm,height=7cm,
keepaspectratio=true]{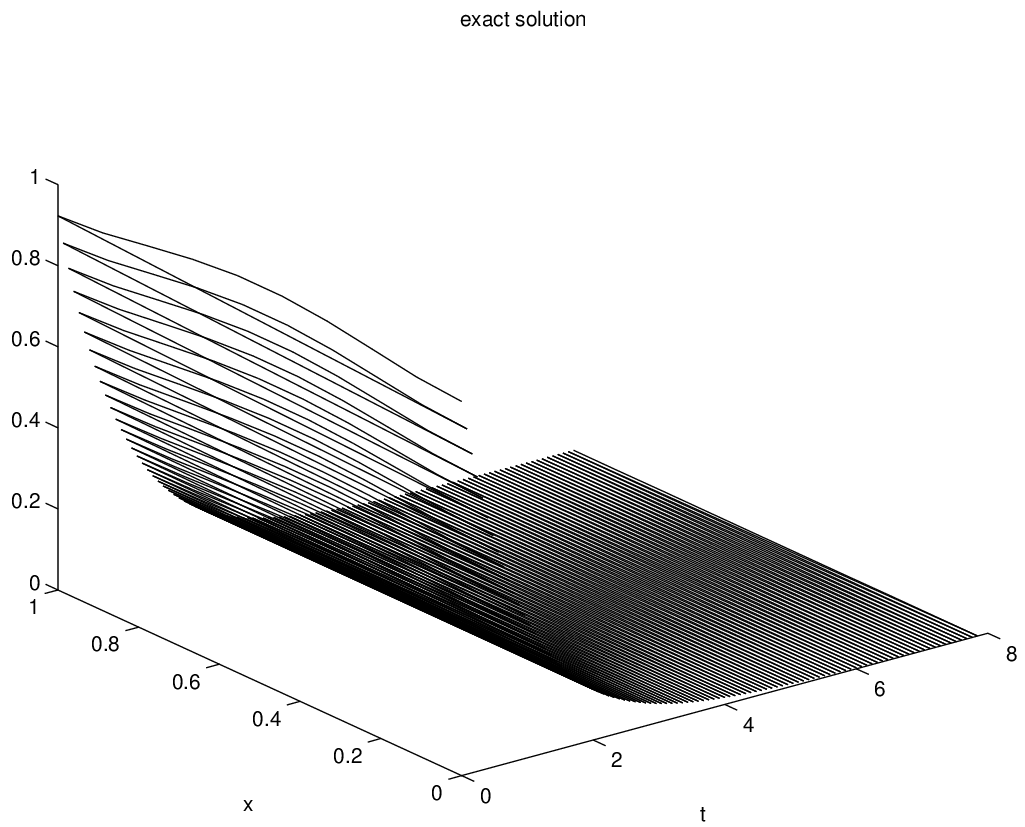}

Figure 2: Exact solution
\end{center}

The fig.3 corresponds to the surface $(x,t)\mapsto u(x,t)$
approximated solution in the case where $f(x,t)=0$. So in both cases
we notice the very good decay of these surfaces from  $T=0$ to
$T=8$. \\

\begin{center}

\includegraphics[width=7cm,height=7cm,
keepaspectratio=true]{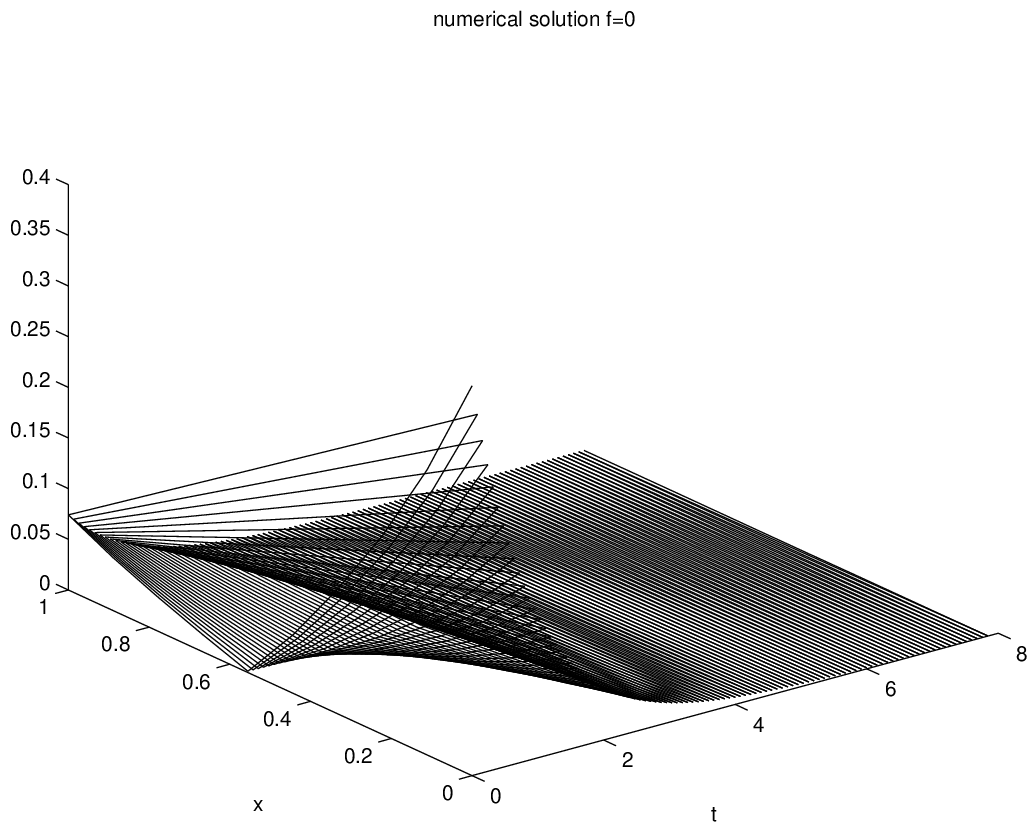}

Figure 3: case f=0
\end{center}


\begin{thebibliography}{}
\bibitem{Clark} H. Clark, {\em Global classical solutions to the Cauchy problem for a nonlinear wave equation},
Internat. J. Math. and Math. Sci. {\bf {21}}(3), (1998), 533-548.
\bibitem{Medeiros} L. A. Medeiros, J. Limaco, S. B. Menezes, {\em Vibrations of Elastic String: Mathematical Aspects, Part one}, J. Comput.  Anal. Appl. {\bf {4}}(2), 2002, 91-127.
\bibitem{Medeiros} L. A. Medeiros, J. Limaco, S. B. Menezes, {\em Vibrations of Elastic String: Mathematical Aspects, Part two}, J. Comput.  Anal. Appl. {\bf {4}}(3) , 2002, 211-263.
\bibitem{Menzala} G. P. Menzala, {\em On global classical solutions of a nonlinear wave equation}, Appl. Anal.
{\bf {10}}, (1980), 179-195.
\bibitem{Lions} J. L. Lions, {\em Quelques m\'{e}thodes de r\'{e}solution des probl\`{e}mes aux limites nonlin\'{e}aires}, Dunod; Gauthier-Villars (Paris),1969.
\bibitem{Long1} N. T. Long, A. P. N. Dinh, {\em On the quasilinear wave equation $u_{tt}- \Delta u + f(u, u_{t}) = 0$ associated with a mixed nonhomogeneous condition}, Nonlinear Anal. {\bf{19}} (1992), 613-623.
\bibitem{Long3} N. T. Long, V. G. Giai, {\em  Existence and asymptotic expansion for a nonlinear wave equation associated with nonlinear boundary conditions}, Nonlinear Anal. Series A: Theory and Methods, {\bf{67}} (6) (2007), 1791-1819.
\bibitem{Santos} M. L. Santos, {\em Decay rates for solutions of a system of wave equations with memory}, EJDE,  Vol {\bf 2002} (2002), No. 38, p. 1 -17.
\end{thebibliography}
\end{document}